\documentstyle[amsfonts,11pt]{article}

 \makeatother

\newtheorem{thm}{Theorem}[section]

\begin{document}
\title{The $\omega$-Lie algebra defined by the commutator of an $\omega$-left-symmetric algebra is not perfect}

\author{Zhiqi Chen  \\ School of Mathematics and Statistics, \\ Guangdong University of Technology, Guangzhou 510520, P.R. China. \\ E-mail: chenzhiqi@nankai.edu.cn \\ Junna Ni \\ Corresponding author. Department of Mathematics, \\ South China Normal University, Guangzhou 510631, P.R.China. \\ Email: nijunna@126.com \\ Jianhua Yu \\ Department of Mathematics, South China Normal University, \\ Guangzhou 510631, P.R.China. Email: yujianhuscnu@126.com}

\maketitle

\begin{abstract}
In this paper, we study admissible $\omega$-left-symmetric algebraic structures on $\omega$-Lie algebras over the complex numbers field $\mathbb C$. Based on the classification of $\omega$-Lie algebras, we prove that any perfect $\omega$-Lie algebra can't be the $\omega$-Lie algebra defined by the commutator of an $\omega$-left-symmetric algebra.


\end{abstract}

\section{Introduction}
A vector space $L$ over $\mathbb F$ is called an $\omega$-Lie algebra if there is a bilinear map $[\cdot,\cdot]:L\times L\rightarrow L$ and a
skew-symmetric bilinear form $\omega: L \times L \rightarrow \mathbb F$ such that
\begin{enumerate}
  \item $[x,y]=-[y,x]$,
  \item $[[x,y],z]+[[y,z],x]+[[z,x],y]=\omega(x,y)z+\omega(y,z)x+\omega(z,x)y$,
\end{enumerate}
hold for any $x,y,z\in L$. Clearly, an $\omega$-Lie algebra $L$ with $\omega=0$ is a Lie algebra, which is called a trivial $\omega$-Lie algebra. Otherwise, $L$ is called a nontrivial $\omega$-Lie algebra. Similar to Lie algebraic case, an $\omega$-Lie algebra $L$ is called perfect if $[L,L]=L$. The notation of an $\omega$-Lie algebra is given by Nurowski in \cite{Nur}, and then there are many studies in this field such as \cite {chen1, 4dim,  chen3,cny, zhang2, Zu}.

For an $\omega$-Lie algebra $L$, a finite dimensional vector space $V$ is called an $L$-module if there exists a bilinear map $L\times V\rightarrow V$ defined by $(x,v)\mapsto xv$ such that
$$[x,y]v=x(yv)-y(xv)+\omega(x,y)v$$
holds for $x,y\in L$ and $v\in V$. It is well-known that left-symmetric algebras are defined by the module of Lie algebras. Similarly, $\omega$-left-symmetric algebras are defined (\cite{zhang}) as follows. Let $L$ be a vector space over $\mathbb F$ with a bilinear map $(x, y) \mapsto xy$. If there is a
bilinear map $\omega:L\times L\rightarrow \mathbb F$ such that
\begin{equation}
(xy)z-x(yz)-(yx)z + y(xz)=\omega(x,y)z,\quad\forall x, y, z\in L,
\end{equation}
then $L$ is called an $\omega$-left-symmetric algebra. Clearly an $\omega$-left-symmetric algebra $L$ with $\omega=0$ is a left-symmetric algebra, which is called a trivial $\omega$-left-symmetric algebra. Otherwise, $L$ is called a nontrivial $\omega$-left-symmetric algebra. Left-symmetric algebras (or pre-Lie algebras, quasi-associative algebras, Vinberg algebras and so on) are first introduced by A. Cayley in 1896 (\cite{Cay}). They appear in many fields in mathematics and mathematical physics, for more details see \cite{Ba, Bor, Bur,Cha, Chu, Eb, Et, Ge, Me, Vin} and so on.

For an $\omega$-left-symmetric algebra $L$, define the commutator $$[x,y]=xy-yx,$$ then $L$ is an $\omega$-Lie algebra. Moreover the left multiplication $l_x$ on the $\omega$-left-symmetric algebra $L$ defined by $l_xy=xy$ makes $L$ to be an $\omega$-Lie algebra $L$-module.
It is an important result that the Lie algebra defined by the commutator of a left-symmetric algebra is not a perfect Lie algebra (\cite{He1}).
In this paper, we prove this result holds for $\omega$-Lie algebras and $\omega$-left-symmetric algebras, i.e.
\begin{thm}\label{lsa}
Let $L$ be an $\omega$-left-symmetric algebra. Then the $\omega$-Lie algebra defined by the commutator $[x,y]=xy-yx$ is not a perfect $\omega$-Lie algebra.
\end{thm}

Throughout this paper, all vector spaces and algebras are finite dimensional over $\mathbb C$ unless stated otherwise.

\section{Perfect $\omega$-Lie algebras}

By the classification of nontrivial $\omega$-Lie algebras, a nontrival perfect $\omega$-Lie algebra $L$ is one of the following cases:
\begin{enumerate}
  \item $\dim L=3$, and $L$ is $A_{\alpha}$, or $B$, or $C_{\alpha}$ for $\alpha\not=0,-1$ in \cite{chen1}. There exists a basis $\{x,y,z\}$ of $L$ such that the nonzero brackets and $\omega$ are given as follows. \\   
   $A_{\alpha}: [x, y] = x, [x,z] = x + y, [y,z] = z + \alpha x, \omega(y,z) = -1$, where $\alpha \in \mathbb C$. \\   
   $B: [x, y] = y, [x,z] = y + z, [y,z] = x, \omega(y,z) =2.$ \\   
   $C_{\alpha}: [x, y] = y, [x,z] = \alpha z, [y,z] = x, \omega(y,z) = 1 + \alpha$, where $0, -1\not=\alpha\in \mathbb C$.
   \item $\dim L=4$, and $L$ is $G_{1,\alpha}$, or $H_{1,\alpha}$, or $\widetilde{A_{\alpha}}$, or $\widetilde{B}$, or $\widetilde{C_{\alpha}}$ for $\alpha\not=0,-1$ in \cite{4dim}. there exists a basis $\{x,y,z,e\}$ of $L$ such that the nonzero brackets and $\omega$ are \\
       $G_{1,\alpha}: [e, x]=e+\alpha y, [e, y]=-e+x, [y,z] = z, [x, y] = y, \omega(e, x)=\alpha, \omega(x, y) = 1.$ \\
       $H_{1,\alpha}: [e, x] = e + \alpha y,[e, y]=-e + x + z, [y,z] = z, [x, y] = y,  \omega(e, x) = \alpha, \omega(x, y) = 1.$ \\
       $\widetilde{A_{\alpha}}: [x,y]=x, [x,z]=x+y, [y,z]=z+{\alpha}x,  [e,z]=e, \omega(y,z)=-1$.\\
       $\widetilde{B}: [x,y]=y, [x,z]=y+z, [y,z]=x, [e,x]=-2e, [e,y]=-e, \omega(y,z)=2.$\\
       $\widetilde{C_{\alpha}}\ ({\alpha}\neq 0,-1): [x,y]=y, [x,z]={\alpha}z, [y,z]=x, [e,x]=-(1+{\alpha})e, \omega(y,z)=1+{\alpha}.$
  \item $\dim L\geq 5$, and $L={\mathbb C}h_0\oplus H_1\oplus {\mathbb C}x\oplus {\mathbb C}v$, and the nonzero brackets and $\omega$ are given as follows: for any $h\in H_1$,
$$[x,h_0]=-ah_0, \ [v,h]=\frac{1}{a}h,\ [v,h_0]=h_2+\frac{1}{a}h_0+x, \ [x,v]=h_1+av,\ \omega(x,v)=1,$$ where $h_1\in {\mathbb C}h_0\oplus H_1$, $h_2\in H_1$ and $a\not=0$. It is type-$P_1$ in \cite{cny}.
  \item $\dim L\geq 5$, and $L=H\oplus {\mathbb C}x\oplus {\mathbb C}y\oplus{\mathbb C}a$, and the nonzero brackets and $\omega$ are given as follows: for any $h\in H$,
$$[a,h]=h, \ [x,y]=h_3+a,\ [x,a]=h_1+b_1x+b_2y,\ [y,a]=h_2+c_1y,\ \omega(x,y)=1,$$
where $h_1,h_2,h_3\in H$, $b_1\not=0$, $c_1\not=0$ and $b_1+c_1+1=0$. It is type-$P_2$ in \cite{cny}.
\end{enumerate}

\section{$\omega$-left-symmetric algebraic structure on perfect $\omega$-Lie algebras.}

This section is to prove Theorem~\ref{lsa}. That is, we extend the classical result to $\omega$-Lie algebras: a Lie algebra admitting a left-symmetric algebraic structure is not a perfect Lie algebra.

Assume that $L$ is an $\omega$-left-symmetric algebra such that the $\omega$-Lie algebra defined by the commutator $[x,y]=xy-yx$ is a perfect $\omega$-Lie algebra. Denote it also by $L$. It is enough to discuss the case when $L$ is not a Lie algebra. Thus $L$ is someone in Section 2. Define $l_u: L \rightarrow L$ by $l_u(v)=uv.$ By the definition of $\omega$-left-symmetric algebra,
 $$l_{[u,v]}=[l_u,l_v]+\omega(u,v)id,\quad \forall u,v\in L.$$
We will discuss perfect $\omega$-Lie algebras in Section 2. case by case.

{\bf Case 1: $\dim L=3$.} In \cite{Nur}, Nurowski classified nontrivial $\omega$-Lie algebras in dimension 3 over $\mathbb R$. Based on the classification, the classification of $\omega$-left-symmetric algebras in dimension 3 is given in \cite{chen} as follows: if $L$ is a nontrivial $\omega$-left-symmetric algebra over ${\mathbb R}$ in dimension 3, then $V$ has a basis $\{e_1,e_2,e_3\}$ such that one of the following cases holds:
\begin{enumerate}
  \item $e_1e_1=e_2e_1=a_1e_1+a_2e_2+a_3e_3=2e_1-e_3e_1$, $e_1e_2=e_2e_2=(a_1-1)e_1+(a_2+1)e_2+a_3e_3=2e_3-e_3e_2$, $e_1e_3=e_2e_3=(2-a_1)e_1+(1-a_2)e_2+(1-a_3)e_3=2e_3-e_3e_3$, $\omega(e_1,e_2)=0$, $\omega(e_2,e_3)=2$, $\omega(e_3,e_1)=-2$.
  \item $e_1e_1=2e_1,e_1e_2=2e_2, e_1e_3=2e_3$, $e_2e_1=e_3e_1=e_2+e_3$, $e_2e_2=e_3e_2=a_1e_1+a_2e_2+a_3e_3$, $e_2e_3=e_3e_3=(a_1+1)e_1+a_2e_2+a_3e_3$, $\omega(e_1,e_2)=0$, $\omega(e_2,e_3)=2$, $\omega(e_3, e_1)=0$.
\end{enumerate}
Then by the above classification of $\omega$-left-symmetric algebras in dimension 3 over $\mathbb R$, we only need to discuss $L=A_{\alpha}$ and $L=C_{\alpha}$.

For the case $L=A_{\alpha}$, the nonzero brackets and $\omega$ are given as follows:
$$[x,y]=x,\quad [x,z]=x+y,\quad [y,z]=z+\alpha x, \quad \omega(y,z)=-1.$$
Then we have:
$$l_x=[l_x,l_y],\quad l_x+l_y=[l_x,l_z],\quad l_z+\alpha l_x=[l_y,l_z]-id.$$
It follows that
\begin{eqnarray*}
    l_x+l_y & = & [l_x,l_z]=[[l_x,l_y],l_z]=[[l_x,l_z],l_y]+[l_x,[l_y,l_z]]  \\
            & = & [l_x,l_y]+[l_x,l_z] \\
            & = & 2l_x+l_y.
\end{eqnarray*}
That is, $l_x=0$. It follows that $l_y=0$ and $l_z=-id$. Then $x=[x,y]=l_xy-l_yx=0$, which is impossible.

For the case $L=C_{\alpha}$, the nonzero brackets and $\omega$ are given as follows:
$$[x,y]=y,\quad [x,z]=\alpha z,\quad [y,z]=x, \quad \omega(y,z)=1+\alpha,$$
where $\alpha\not=0,-1$. Then we have:
$$l_y=[l_x,l_y],\quad \alpha l_z=[l_x,l_z],\quad l_x=[l_y,l_z]+(1+\alpha)id.$$
It follows that
\begin{eqnarray*}
    l_x & = & [l_y,l_z]+(1+\alpha)id=[[l_x,l_y],l_z]+(1+\alpha)id   \\
            & = & [[l_x,l_z],l_y]+[l_x,[l_y,l_z]]+(1+\alpha)id \\
            & = & \alpha[l_z,l_y]+(1+\alpha)id \\
            & = & -\alpha l_x+(1+\alpha)^2id.
\end{eqnarray*}
That is, $l_x=(1+\alpha)id$ since $\alpha\not=-1$. It follows that $l_y=l_z=0$. Then $x=[y,z]=l_yz-l_zy=0$, which is impossible.

{\bf Case 2: $\dim L=4$.} we will discuss case by case. {\bf The first case is $L=G_{1,\alpha}$.} Then the nonzero brackets and $\omega$ are given as follows:
\begin{eqnarray*}
&& [e,x]=e+\alpha y, \quad [e,y]=-e+x, \quad [y,z]=z, \quad [x,y]=y, \\
&& \omega(e,x)=\alpha, \quad \omega(x,y)=1.
\end{eqnarray*}
Let $e'=e+\alpha y.$ Then for $\{e',x,y,z\}$, we have:
\begin{eqnarray*}
   && [e',x]=e'-\alpha y, \quad [e',y]=-e'+x+\alpha y, \quad  [e',z]={\alpha}z, \\
   && [y,z]=z, \quad [x,y]=y, \quad \omega(x,y)=1.
\end{eqnarray*}
Then we have:
\begin{eqnarray*}
   && l_y=[l_x,l_y]+id, \quad [l_x,l_z]=0, \quad l_z=[l_y,l_z],\\
   && l_{e'}-{\alpha}l_y=[l_{e'},l_x], \quad -l_{e'}+l_x+{\alpha}l_y=[l_{e'},l_y], \quad {\alpha}l_z=[l_{e'},l_z].
\end{eqnarray*}
Clearly $-l_{e'}+l_x+{\alpha}l_y=[l_{e'},l_y]$ means $l_{e'}=l_x+{\alpha}l_y-[l_{e'},l_y]$. Then
\begin{eqnarray*}
   l_{e'}-{\alpha}l_y & = &  [l_x+al_y-[l_{e'},l_y],l_x] \\
                      & = &   {\alpha}[l_y,l_x]-[[l_{e'},l_y],l_x]
\end{eqnarray*}
It is easy to see that
\begin{eqnarray*}
   [[l_{e'},l_y],l_x] & = & [[l_{e'},l_x],l_y]+[l_{e'},[l_y,l_x]] = [l_{e'},l_y]-[l_{e'},l_y] \\
                      & = & 0.
\end{eqnarray*}
That is, $l_{e'}-al_y={\alpha}[l_y,l_x]={\alpha}(id-l_y).$ It means
  $$l_{e'}={\alpha}id .$$
Then we have $${\alpha}l_y=aid, \quad {\alpha}l_z=0, \quad l_x=0.$$
Then $e'-{\alpha}y=[e',x]=e'x-xe'=l_{e'}x-l_xe'={\alpha}x$, which is impossible.

{\bf The second case is $L=H_{1,{\alpha}}$.} For this case, the nonzero brackets and $\omega$ are
\begin{eqnarray*}
   && [x,y]=y, \quad [y,z]=z, \quad [e,y]=x+z-e, \quad [e,x]=e+{\alpha}y, \\
   && \omega(x,y)=1, \quad \omega(e,x)={\alpha}.
\end{eqnarray*}
It follows that
\begin{eqnarray*}
   && l_y=[l_x,l_y]+id, \quad l_z=[l_y,l_z], \quad  [l_x,l_z]=0, \\
   && l_x+l_z-l_e=[l_e,l_y], \quad l_e+{\alpha}l_y=[l_e,l_x]+{\alpha}id, \quad [l_e,l_z]=0.
\end{eqnarray*}
Clearly,
\begin{eqnarray*}
   0 & = & [l_x,l_z]=[l_x,[l_y,l_z]]=[[l_x,l_y],l_z]+[l_y,[l_x,l_z]] \\
     & = & [l_y,l_z]=l_z.
\end{eqnarray*}
Clearly, $l_e=[l_e,l_x]+{\alpha}id-{\alpha}l_y.$ Then we have
\begin{eqnarray*}
     l_x-l_e&=&[l_e,l_y]=[[l_e,l_x]+{\alpha}id-{\alpha}l_y,l_y]=[[l_e,l_x],l_y]\\
            & =& [[l_e,l_y],l_x]+[l_e,[l_x,l_y]] \\
            & = & -[l_e,l_x]+[l_e,l_y].
\end{eqnarray*}
It gives $l_e={\alpha}(id-l_y).$ Then $$l_x=l_e+[l_e,l_y]=l_e.$$
Furthermore, $l_y-id=[l_x,l_y]=[l_e,l_y]=0$. It follows that $$l_y=id,\quad l_x=l_e=0.$$
Then $y=[x,y]=l_xy-l_yx=-x$, which is impossible.

{\bf The third case is $L=\widetilde{A_{\alpha}}$.} The nonzero brackets and $\omega$ are given
$$[x,y]=x,\quad [x,z]=x+y,\quad [y,z]=z+{\alpha}x, \quad [e,z]=e, \quad \omega(y,z)=-1.$$
Then we have
\begin{eqnarray*}
   && l_x=[l_x,l_y], \quad l_x+l_y=[l_x,l_z], \quad l_z+{\alpha}l_x=[l_y,l_z]-id, \\
   && l_e=[l_e,l_z], \quad [l_e,l_x]=0, \quad [l_e,l_y]=0.
\end{eqnarray*}
Clearly $l_x=[l_x,l_z]-l_y.$ Then
\begin{eqnarray*}
    l_x &=&  [l_x,l_y]=[[l_x,l_z]-l_y,l_y]=[[l_x,l_z],l_y] \\
        &=& [[l_x,l_y],l_z]+[l_x,[l_z,l_y]] \\
        &=& [l_x,l_z]+[l_x,-l_z-{\alpha}l_x-id]  \\
        &=& 0.
\end{eqnarray*}
Then $l_y=0$ by $l_x+l_y=[l_x,l_z]$. But $x=[x,y]=xy-yx=l_x(y)-l_y(x)=0$, which is impossible.

{\bf The fourth case is $L=\widetilde{B}$.} The nonzero brackets and $\omega$ are
$$[x,y]=y,\quad [x,z]=y+z, \quad [y,z]=x, \quad [e,x]=-2e, \quad [e,y]=-e,\quad  \omega(y,z)=2.$$
Then we have
\begin{eqnarray*}
   &&  l_y=[l_x,l_y], \quad  l_y+l_z=[l_x,l_z],\quad  l_x=[l_y,l_z]+2id, \\
   && -2l_e=[l_e,l_x], \quad -l_e=[l_e,l_y], \quad [l_e,l_z]=0.
\end{eqnarray*}
Since $l_e=[l_y,l_e]$, we have
\begin{eqnarray*}
-2l_e &=& [l_e,l_x]=[[l_y,l_e],l_x]=[[l_y,l_x],l_e]+[l_y,[l_e,l_x]] \\
      &=& [-l_y,l_e]+[l_y,-2l_e]=-3[l_y,l_e] \\
      &=& -3l_e.
\end{eqnarray*}
That is $l_e=0.$ Since $l_y=[l_x,l_y]$, we have
\begin{eqnarray*}
     l_x &=& [l_y,l_z]+2id=[[l_x,l_y],l_z]+2id=[[l_x,l_z],l_y]+[l_x,[l_y,l_z]]+2id \\
         &=& [l_y+l_z,l_y]+[l_x,l_x-2id]+2id  \\
         &=& 4id-l_x.
\end{eqnarray*}
It means $l_x=2id$. Then $l_y=l_z=0$. But $y=[x,y]=l_xy-l_yx=2y$, which is impossible.

{\bf The fifth case is $L=\widetilde{C_{\alpha}}\ ({\alpha}\neq 0,-1)$.} The brackets and $\omega$ are given as follows:
$$[x,y]=y, \quad [x,z]={\alpha}z, \quad [y,z]=x, \quad [e,x]=-(1+{\alpha})e, \quad \omega(y,z)=1+{\alpha}.$$
Then we have
\begin{eqnarray*}
    && l_y=[l_x,l_y], \quad {\alpha}l_z=[l_x,l_z], \quad l_x=[l_y,l_z]+(1+{\alpha})id \\
    && ({\alpha}+1)l_e=[l_x,l_e], \quad [l_e,l_y]=0, \quad [l_e,l_z]=0.
\end{eqnarray*}
Furthermore we have
\begin{eqnarray*}
  l_x&=&[l_y,l_z]+(1+{\alpha})id=[[l_x,l_y],l_z]+(1+{\alpha})id \\
     &=&[[l_x,l_z],l_y]+[l_x,[l_y,l_z]]+(1+{\alpha})id  \\
     &=&{\alpha}[l_z,l_y]+(1+{\alpha})id  \\
     &=& -{\alpha}l_x+(1+{\alpha})^2id.
\end{eqnarray*}
That is $l_x=(1+{\alpha})id.$ Then $l_y=0$ by $l_y=[l_x,l_y]$. But $y=[x,y]=l_xy-l_yx=(1+{\alpha})y$, which is impossible since $\alpha\not=0$.

{\bf Case 3: $P_1$-type.} 
Let $x'=x+h_2$. Then $[v,h_0]=\frac{1}{a}h_0+x'$, $[x',h_0]=-ah_0$ and
$[x',v]=h_1'+av$. Here $h_1'=h_1-\frac{1}{a}h_2.$ Then we have
\begin{eqnarray*}
   &&\frac{1}{a}l_{h_0}+l_{x'}=[l_v,l_{h_0}], \quad l_{h_1'}+al_v+id=[l_{x'},l_v], \quad  -al_{h_0}=[l_{x'},l_{h_0}].
\end{eqnarray*}
Then we have
\begin{eqnarray*}
 [l_v,-al_{h_0}]&=&[l_v,[l_{x'},l_{h_0}]]=[[l_v,l_{x'}],l_{h_0}]+[l_{x'},[l_v,l_{h_0}]]\\
                   &=&[-l_{h_1'}-al_v-id,l_{h_0}]+[l_{x'},\frac{1}{a}l_{h_0}+l_{x'}]\\
                   &=&-a[l_v,l_{h_0}]+\frac{1}{a}[l_{x'},l_{h_0}].
\end{eqnarray*}
It follows that $\frac{1}{a}[l_{x'},l_{h_0}]=-l_{h_0'}=0$, and then $l_{x'}=0.$ Then $-ah_0=[x',h_0]=l_{x'}(h_0)-l_{h_0}(x')=0,$
which is impossible.

{\bf Case 4: $P_2$-type.} For this case, we have
\begin{eqnarray*}
 && [l_x,l_h]=0, \quad [l_y,l_h]=0, \quad l_h=[l_a,l_h], \quad l_{h_1}+b_1l_x+b_2l_y=[l_x,l_a], \\
 && l_{h_2}+c_1l_y=[l_y,l_a],\quad l_{h_3}+l_a=[l_x,l_y]+id.
\end{eqnarray*}
By $l_{h_2}+c_1l_y=[l_y,l_a]$ and $[l_x,l_h]=0,$ we have
\begin{eqnarray*}
[l_x,c_1l_y] & = & [l_x,[l_y,l_a]-l_{h_2}]=[l_x,[l_y,l_a]]-[l_x,l_{h_2}]=[l_x,[l_y,l_a]]-[l_x,l_{h_2}] \\
            & = & [[l_x,l_y],l_a]+[l_y,[l_x,l_a]]
\end{eqnarray*}
Then by $l_{h_3}+l_a=[l_x,l_y]+id $ and $l_{h_1}+b_1l_x+b_2l_y=[l_x,l_a]$, we have
\begin{eqnarray*}
[l_x,c_1l_y] & = & [l_{h_3}+l_a-id, l_a]+[l_y,l_{h_1}+b_1l_x+b_2l_y] \\
            & = & -l_{h_3}-b_1[l_x,l_y].
\end{eqnarray*}
Since $c_1+b_1+1=0$, we have that $[l_x,l_y]=l_{h_3}$.
By $l_{h_3}+l_a=[l_x,l_y]+id$, $$l_{a}=id.$$
It follows that $$l_h=l_y=l_x=0.$$
Then $h_3+a=[x,y]=xy-yx=l_xy-l_yx=0,$ which is impossible.

{\bf Proof of Theorem~\ref{lsa}.} In summary, we finish Theorem~\ref{lsa}.

\section*{Acknowledgements}
Z. Chen was partially supported by NNSF of China (11931009 and 12131012).

\end{document}